\documentclass[12pt,letterpaper]{amsart}

\allowdisplaybreaks
\usepackage{mathpazo}
\usepackage{amsmath,amsfonts,amssymb}
\usepackage{amsrefs}
\usepackage{amsthm}
\usepackage{latexsym,amsmath,amssymb,amsfonts}
\usepackage{rotating}
\usepackage{mathrsfs}
\usepackage{xypic} \xyoption{all}
\usepackage{amscd}
\usepackage{hyperref} 
\usepackage{euscript}
\usepackage{hhline}
\usepackage{graphicx,epstopdf}
\usepackage{epsfig}
\usepackage{xcolor}
\usepackage{textcomp}
\usepackage[all,color]{xy}

\newlength{\fighskip} \fighskip=2pt
\newlength{\figvskip} \figvskip=3pt

\usepackage{hyperref}

\advance\textheight by \topskip
\numberwithin{equation}{section}

\newcommand{\C}{\mathbb{C}}

\newcommand{\W}{\mathcal{W}}

\DeclareMathOperator{\Sym}{\text{Sym}}
\DeclareMathOperator{\Hom}{\text{Hom}}

\DeclareMathOperator{\Tr}{\text{Tr}}

\newcommand{\A}{\mathcal A}

\theoremstyle{plain}
\newtheorem{thm}{Theorem}[section]
\newtheorem{thm-defn}{Theorem/Definition}[section]
\newtheorem{lem}[thm]{Lemma}
\newtheorem{lem-defn}[thm]{Lemma/Definition}
\newtheorem{prop}[thm]{Proposition}
\newtheorem{cor}[thm]{Corollary}

\newtheorem{que}[thm]{Question}

\theoremstyle{definition}
\newtheorem{defn}[thm]{Definition}
\newtheorem{notn}[thm]{Notation}

\theoremstyle{remark}
\newtheorem{rmk}[thm]{Remark}

\allowdisplaybreaks[4]  

\begin{document}
\title{Quantization of K\"ahler manifolds}

\author[Chan]{Kwokwai Chan}
\address{Department of Mathematics\\ The Chinese University of Hong Kong\\ Shatin\\ Hong Kong}
\email{kwchan@math.cuhk.edu.hk}

\author[Leung]{Naichung Conan Leung}
\address{The Institute of Mathematical Sciences and Department of Mathematics\\ The Chinese University of Hong Kong\\ Shatin \\ Hong Kong}
\email{leung@math.cuhk.edu.hk}

\author[Li]{Qin Li}
\address{Department of Mathematics\\ Southern University of Science and Technology\\ Shenzhen\\China}
\email{liqin@sustech.edu.cn}


\dedicatory{Dedicated to Professor Mathai Varghese on the occasion of his 60th birthday.}

\begin{abstract}
	This is a survey on our recent works \cites{Chan-Leung-Li, CLL, CLL3} which reveal new relationships among deformation quantization, geometric quantization, Berezin-Toeplitz quantization and BV quantization on K\"ahler manifolds.
\end{abstract}

\maketitle

\section{Introduction}

Quantization is an important subject in symplectic geometry, mathematical physics, as well as in representation theory. There are various quantization schemes including deformation quantization, geometric quantization, Berezin-Toeplitz quantization and quantum field theories such as Batalin-Vilkovisky (BV) quantization. In this survey article, we will explain our recent works \cites{Chan-Leung-Li, CLL, CLL3} on new relationships among them. We will concentrate on K\"ahler manifolds $( X,\omega ,J) $, namely, symplectic manifolds endowed with complex polarizations.

{\em Geometric quantization} is about producing a Hilbert space $H$ associated to $X$, depending on the quantum parameter $\hbar=1/k$ for $k$ a large integer \cite{Donaldson}. We need to assume that $\left[ \omega \right] $ is integral (i.e., $[ \omega ] \in H^{2}( X,\mathbb{Z}) $). The procedure would depend on a choice of {\em polarization}, but the resulting $H$ should be canonical up to scalings. When $( X,\omega) $ has a Hamiltonian symmetry group $G$ and there is an invariant polarization, the group $G$ would act on the corresponding $H$. Conversely, Kirillov's program says that most important representations of $G$ should arise in this way.

On the other hand, {\em deformation quantization} is a formal deformation of the commutative algebra $(C^{\infty }(X),\cdot)$ equipped with pointwise multiplication to a noncommutative one $( C^{\infty }( X)[[\hbar]] ,\star_{\hbar}) $ equipped with a {\em star product}, so that the leading order of noncommutativity is given by the Poisson bracket $\{-,-\} $, i.e.,
\begin{equation}\label{equation: Poisson-bracket}
\frac{d}{d\hbar}\Big|_{\hbar=0}\left( f\star_{\hbar}g - g\star_{\hbar}f\right)
=\left\{ f,g\right\}.
\end{equation}
Here $\left\{ f,g\right\} $ is the Lie algebra structure on $C^{\infty}\left( X\right) $ induced from the Lie algebra $\text{Vect}( X,\omega)$ of symplectic vector fields via $f\mapsto X_{f}$, where $X_{f}$ is the Hamiltonian vector field determined by $\omega \left( X_{f},-\right) =df$.

When $X=T^*M$ equipped with the canonical symplectic form $\omega=\sum_j dx^{j}\wedge dp_{j}$ where $x^{j}$'s are local coordinates of $M$ and $p_{j}$'s are corresponding cotangent coordinates, we could choose $H=L^{2}(M)$, the space of $L^{2}$-functions on $M$. To view $s(x) \in H$ globally on $X$, they are certain functions on $X$ which are constant along the $p$-directions. The span $\mathcal{D}=\left\langle \partial/\partial p_1, \dots, \partial/\partial p_n \right\rangle \subset TX$, where $n = \dim M$, is an integrable Lagrangian distribution of $X$, which is a {\em real polarization}. Since $X=T^*M$, the trivial line bundle $L$ over $X$ admits a Hermitian metric with a unitary connection $D_{A}$ whose curvature is $F_{A}=-2\pi\sqrt{-1}\omega $. Therefore $H$ can be identified as the space $\Gamma _{\mathcal{D}}\left( X,L\right) $ of $L^{2}$-sections of $L$ over $X$ which are covariantly constant along $\mathcal{D}$.

To explain deformation quantization, we consider $X=T^*\mathbb{R}^{n}$. Given $f(x,p) \in C^{\infty}( T^*\mathbb{R}^{n})$, if its dependency on $p_{j}$'s is polynomial, i.e., $f \in \Gamma( \mathbb{R}^{n},\Sym^* T_{\mathbb{R}^{n}})$, then we could view it as the symbol of a linear differential operator $\hat{f}$
operating on functions $s(x) $ by treating $p_{j}$ as $\hat{p}_{j} = \sqrt{-1}\hbar\frac{\partial }{\partial x^{j}}$. The deformed product $f \star_{\hbar} g$ becomes the composition of the differential operators $\hat{f}$ and $\hat{g}$ and equation \eqref{equation: Poisson-bracket} is just the uncertainty principle%
\[
[ \hat{x}^{j},\hat{p}_{k}] \sim \sqrt{-1}\hbar\delta _{k}^{j}.
\]

On a general symplectic manifold $X$, equipped with a real polarization $\mathcal{D}$, we could set $H = \Gamma _{\mathcal{D}}(X, L)$ in the same manner as above provided that such an $L$ exists. $L$ is called a
{\em prequantum line bundle} (i.e., $c_1(L) = [\omega]$) and its existence is guaranteed if $[\omega]$ is integral. It may not be easy to construct real polarizations on a general symplectic manifold. But if we allow Lagrangian foliations $\mathcal{D}_{c}$ on the complexified tangent bundle, i.e., $TX\otimes \mathbb{C}=\mathcal{D}_{c}\oplus \mathcal{\bar{D}}_{c}$, then such structures are equivalent to K\"ahler structures on $X$ with holomorphic tangent bundle given by $T^{1,0}X=\mathcal{D}_{c}$. These are called {\em complex polarizations} on $( X,\omega)$. The integrality condition on $[\omega]$ is the same as $X$ being a complex projective manifold by the
Kodaira embedding theorem. In this case, $H=\Gamma _{\mathcal{D}_{c}}(X, L) = H^{0}(X, L)$, the space of $L^{2}$-holomorphic sections of $L$.

For a K\"ahler manifold $( X,\omega ,J)$, taking holomorphic sections $H_{k}=H^{0}( X,L^{\otimes k}) $ is the recipe of producing geometric quantization when $[ \omega ] $ is integral. For every smooth function $f$ on $X$, there is an associated operator given by the multiplication $m_{f}$ by $f$ followed by the orthogonal projection $\Pi_k :\Gamma _{L^{2}}( X,L^{\otimes k})
\rightarrow H^{0}( X,L^{\otimes k}) $ and this is called the {\em Toeplitz operator} $T_{f,k} := \Pi_k \circ m_{f}$. Then there are bi-differential operators $C_i(-,-)$ and constants $K_N(f,g)$, which are independent of $k$, such that the following estimates hold:
\begin{equation}\label{equation: asymptotic-expansion-composition-Toeplitz-operators}
\left\Vert T_{f,k}\circ T_{g,k} - \sum_{i=0}^{N-1}\left(\frac{1}{k}\right)^iT_{C_i(f,g),k}\right\Vert\leq K_N(f,g)\left(\frac{1}{k}\right)^N.
\end{equation}
Here $\Vert \cdot \Vert$ denotes the operator norm.  This defines the {\em Berezin-Toeplitz star product} $\star_{BT}$ and the {\em Berezin-Toeplitz deformation quantization algebra} $(C^\infty(X)[[\hbar]],\star_{BT})$ \cites{Bordemann-Meinrenken, Schlichenmaier, Karabegov} (see also \cites{Bordemann, Bordemann-Waldmann, DLS, Karabegov96, Karabegov00, Karabegov07, Ma-Ma, Neumaier}). The Berezin-Toeplitz quantization is a special example of Wick type deformation quantization of K\"ahler manifolds.

When $X=\mathbb{C}^{n}$ equipped with the trivial holomorphic line bundle $L$,  we may take $k=1$. Then $H=H^{0}(X, L)$ consists of holomorphic functions $s(z)$ which are $L^{2}$ with respect to $\left\Vert s\right\Vert^{2}=\int_{X}\left\vert s(z) \right\vert ^{2}e^{-\left\vert z\right\vert ^{2}/2}$.  For deformation quantization on $\C^n$, we restrict ourselves to polynomials. Explicitly,  $f( z,\bar{z})$ acts as a linear differential operator $T_f$ acting on $s(z)\in H^0(X,L)$ by letting $T_{{z}_{j}}=z_{j}\cdot $ and $T_{\bar{z}_{j}}=\hbar\frac{\partial }{\partial z_{j}}$ with together with the {\em Wick ordering}. (See section \ref{subsection: flat-space} for more details). Then the composition of these $T_f$'s defines a deformation quantization of polynomials on $\mathbb{C}^{n}$, with $\hbar=1$ not only a formal variable. In this case, the left hand side of the estimate in equation \eqref{equation: asymptotic-expansion-composition-Toeplitz-operators} vanishes if $N>\max\{\deg(f),\deg(g)\}$. This implies that for $X=\C^n$, there is $T_f\circ T_g=T_{f\star_{BT} g}$, and $H$ forms an honest representation of this deformation quantization of polynomials. 


For a general compact K\"ahler manifold $X$, however, the estimate \eqref{equation: asymptotic-expansion-composition-Toeplitz-operators} only says that the difference $T_{f, k} \circ T_{g, k} - T_{f\star_{BT} g, k}$ is {\em asymptotically} zero when $k$ tends to infinity. Thus we only have an ``asymptotic action'' of the Berezin-Toeplitz deformation quantization on the $H_k$'s as $k$ approaches infinity. Although it is enough to recover the star product, none of the Hilbert spaces $H_k$'s actually form a representation of the deformation quantization algebra. It is therefore natural to ask the following question (cf. problem iv in \cite[Section 9]{Bordemann-Waldmann}):

\begin{que}\label{question}
 Do we have an honest Toeplitz action of $C^{\infty}\left(X\right)[[\hbar]] $ on a ``geometric quantization'' so that this would automatically determine a star product $\star_{\hbar}$ as in the flat case? 
\end{que}


In \cites{Chan-Leung-Li, CLL, CLL3} we gave a positive answer to this question: by suitably localizing the Hilbert spaces $H^{0}\left( X,L^{\otimes k}\right)$ using so-called peak sections, we can construct a family of representations $H_{x_0}$ of $C^{\infty}\left( X\right)[[\hbar]]$ parametrized by points in $X$ via a formal Toeplitz action. Combining with Fedosov's approach to deformation quantization, we can patch these $H_{x_0}$'s together to produce a sheaf version of geometric quantization over $X$. Our construction is more natural from a quantum geometric perspective: notice that a star product is local in nature, so a deformation quantization defines the structure sheaf of a ``quantum geometry'' on $X$.
In the rest of this introduction, we will briefly explain our results; more detailed descriptions of our constructions and results will be given in subsequent sections.



Using peak sections technique \cite{Tian}, we localize global holomorphic sections of $L^{\otimes k}$  to approximate the flat model near $x_{0}$. In \cites{Chan-Leung-Li}, we devised a method to link localized sections of $L^{\otimes k}$ with different $k$'s so that we could honestly take the $k\rightarrow \infty $ limit of all these sections to recover the Hilbert space $H_{x_{0}}$ as constructed in the standard flat situation.  As a consequence, we obtained a representation of $C^{\infty}\left( X\right)[[\hbar]] $ on $H_{x_{0}}$ by restricting smooth functions on $X$ to a formal neighborhood of $x_{0}$. More precisely, we proved the following theorem (see Section \ref{section: BT-representation} for more details):

\begin{thm}[Theorem 1.1 in \cites{Chan-Leung-Li}]
For every $x_0\in X$, the vector spaces $H_{x_0}$ constructed above form a representation of the Berezin-Toeplitz deformation quantization $(C^\infty(X)[[\hbar]],\star_{BT})$ satisfying locality, i.e., for every smooth function $f$, its action on $H_{x_0}$ depends only on the infinite jets of $f$ at $x_0$.  
\end{thm}

To put all these $H_{x_{0}}$'s together to form a sheaf over $X$, we mimic Fedosov's approach. In \cites{Fed, Fedbook}, Fedosov gave a geometric construction of deformation quantization 
on symplectic manifolds by gluing the Moyal product on  each tangent space $T_{x}X\simeq T^*\mathbb{R}^{n}\simeq \mathbb{C}^{n}$ using a flat connection. This is now known as a {\em Fedosov abelian connection}. On a K\"ahler manifold $X$, the fiberwise Moyal product is replaced by the polarized Wick product on the complexified Weyl space $\mathcal{W}_{x_{0}}\simeq \mathbb{C}[[z^{1},\bar{z}^{1},\cdots, z^n,\bar{z}^n]][[\hbar]]$. In \cite{CLL}, we constructed a special family of Fedosov abelian connections which are natural quantizations of {\em Kapranov's $L_\infty$ structure} \cite{Kapranov} (see also \cite{Kontsevich2}), and showed that every Wick type star product on $X$ can be obtained by these Fedosov abelian connections (see Section \ref{section: L_infty-structures} for more details):
\begin{thm}[Theorem 1.1 in \cite{CLL}]
Let $\alpha$ be a representative of a formal cohomology class in $\hbar H^2_{dR}(X)[[\hbar]]$ of type $(1,1)$. Then there exists $\gamma_\alpha \in \mathcal{A}_X^{0,1}(\mathcal{W}_{X,\mathbb{C}})$ solving the Fedosov equation
 \begin{equation*}
\nabla \gamma_\alpha + \frac{1}{\hbar} \gamma_\alpha\star \gamma_\alpha + R_\nabla=\alpha,
\end{equation*}
so that $D_{F,\alpha}:=\nabla+\frac{1}{\hbar}[\gamma_\alpha,-]_\star$ defines a Fedosov abelian connection. This is a quantum extension of Kapranov's $L_\infty$ structure $D_K$ in the sense that $D_{F,\alpha}|_{\W_X}=D_K$. Furthermore, the star product associated to $D_{F,\alpha}$ is of Wick type with Karabegov form given by $2\sqrt{-1}\omega-\alpha$. 
\end{thm}

The Wick product on the Weyl space $\mathcal{W}_{x_{0}}\simeq \mathbb{C}[[z^{1},\bar{z}^{1},\cdots, z^n,\bar{z}^n]][[\hbar]]$ acts on the Fock space $H_{x_{0}}\simeq \mathbb{C}[[z^{1},\cdots, z^n]][[\hbar]]$, which is known as the {\em Bargmann-Fock representation}. By gluing these fiberwise representations together using Fedosov's technique, we obtained in \cite{CLL3} a global vector bundle $\mathcal{F}_{X,\alpha }$ over $X$ equipped with a natural flat connection which is compatible with the Fedosov abelian connection on $\mathcal{W}_{X,\C}$ via the above fiberwise action. This compatibility enables us to construct a module sheaf via (local) flat sections:
\begin{thm}[Theorem 1.1 in \cite{CLL3}]
The {\em Bargmann-Fock sheaf} $\mathcal{F}_{X,\alpha}^{\text{flat}}$, which consists of flat sections of $\mathcal{F}_{X,\alpha}$ satisfying a convergence property, is a sheaf of modules over $(C^{\omega}_X[[\hbar]],\star_\alpha)$. 
\end{thm}
Here we need the Toeplitz action of functions  on $\mathcal{F}_{X,\alpha}^{flat}$ to preserve the convergence property, so we restrict to the subspace $C^{\omega}[[\hbar]]$ of formal analytic functions. Furthermore, we showed that this action is indeed given by the formal Toeplitz
action. More precisely, we have the following
\begin{thm}[Theorem 1.2 in \cite{CLL3}]
For every $x_0\in X$, there exists a subspace $V_{x_0}$ of the stalk $\left(\mathcal{F}_{X,\alpha}^{\text{flat}}\right)_{x_0}$ isomorphic to the space of germs of formal holomorphic functions at $x_0$, i.e., $V_{x_0}\cong\mathcal{O}_{X,x_0}[[\hbar]]$, such that the representation of smooth functions on $V_{x_0}$ is via {\em formal Toeplitz operators}. 
\end{thm}
See Section \ref{section: BF-sheaves} for more details.

Now the following corollary gives an answer to Question \ref{question}:
\begin{cor}
Action of the formal Toeplitz operators of $C^{\omega}\left(X\right)[[\hbar]] $ on the
Bargmann-Fock sheaf $\mathcal{F}_{X,\alpha}^{\text{flat}}$ uniquely determines the Wick type star product whose Karabegov form is $2\sqrt{-1}\omega-\alpha$.
\end{cor}

Following \cite{GLL}, we considered in \cite{CLL} the Batalin-Vilkovisky (BV) quantizations of K\"ahler manifolds associated to the Fedosov abelian connections we obtained. By running the homotopy group flow operator to $\gamma_\alpha$, we obtain a solution $\gamma_\infty$ of the quantum master equation (QME). This is the $\infty$-scale effective renormalization of a one-dimensional Chern-Simons theory with target $X$. A significant feature of $\gamma_\infty$ is that, in contrast with the real symplectic case studied in \cite{GLL}, its graph expansion involves only trees and one-loop graphs. So this produces one-loop exact BV quantizations on a K\"ahler manifold $X$. In particular, this gives an explicit, cochain level computation of the trace of Wick type star products on $X$. More details can be found in Section \ref{section: BV}. 

\subsection*{Acknowledgement}
\

We thank Si Li and Siye Wu for useful discussions. The first named author thanks Martin Schlichenmaier and Siye Wu for inviting him to attend the conference GEOQUANT 2019 held in September 2019 in Taiwan, in which he had stimulating and very helpful discussions with both of them as well as Jørgen Ellegaard Andersen, Motohico Mulase, Georgiy Sharygin and Steve Zelditch. The third named author thanks the organizers for inviting him to the conference ``Index Theory, Duality and Related Fields'' held in the Chern Institute of Mathematics in June 2019, where he had very helpful discussions with Xiaonan Ma, Xiang Tang and Siye Wu. 

K. Chan was supported by grants of the Hong Kong Research Grants Council (Project No. CUHK14302617 \& CUHK14303019) and direct grants from CUHK.
N. C. Leung was supported by grants of the Hong Kong Research Grants Council (Project No. CUHK14301117 \& CUHK14303518) and direct grants from CUHK.
Q. Li was supported by grants from National Science Foundation of China (Project No. 12071204), and Guangdong Basic and Applied Basic Research Foundation (Project No. 2020A1515011220).

\section{Berezin-Toeplitz quantization and its representations}\label{section: BT-representation}
Deformation quantization is the mathematical description of quantum observables on phase spaces. From the canonical quantization point of view, these quantum observables are operators on Hilbert spaces. We would like to understand those operators on Hilbert spaces corresponding to smooth functions, or equivalently Hilbert representations of the deformation quantization algebras. 

We focus on phase spaces which are K\"ahler manifolds and Wick type deformation quantizations. 

\subsection{Flat space $\mathbb{C}^n$}\label{subsection: flat-space}
\

First we recall the Toeplitz operators on $X=\C^n$ associated to polynomials. On the flat spaces $ \mathbb{C}^n$ with trivial prequantum line bundle $L$, the Hilbert space on which the Toeplitz operators act is the well-known {\em Bargmann-Fock space} $\mathcal{H}L^2(\mathbb{C}^n,\mu_\hbar)$ consisting of $L^2$ integrable entire holomorphic functions with respect to the density $\mu_\hbar(z) = (\pi\hbar)^{-n}e^{-|z|^2/\hbar}$; here $\hbar$ is regarded as a positive real number.
It is easy to see, by direct computations, that the holomorphic polynomials 
$$
\frac{z^I}{\sqrt{I!\hbar^{|I|}}},
$$
where $I$ runs over all multi-indices, form an orthonormal basis of $\mathcal{H}L^2(\mathbb{C}^n,\mu_\hbar)$. For a monomial $f(z,\bar{z})=z^{i_1}\cdots z^{i_k}\bar{z}^{j_1}\cdots\bar{z}^{j_l}$, the associated Toeplitz operator is given by the following {\em Wick ordering}:
\begin{align*}
T_{f_1(z)f_2(\bar{z})}=\left(\hbar\frac{\partial}{\partial z^{j_1}}\right)\circ\cdots\circ\left(\hbar\frac{\partial}{\partial z^{j_1}}\right)\circ m_{z^{i_1}\cdots z^{i_k}},
\end{align*}
Let $f,g\in\mathbb{C}[z,\bar{z}]$. Then there is a formula for the composition of Toeplitz operators: $T_f\circ T_g=T_{f\star g}$, where 
$$
f\star g:=\exp\left(-\hbar\sum_{i=1}^n\frac{\partial}{\partial z_i}\frac{\partial}{\partial \bar{w}_i}\right)(f(z,\bar{z})g(w,\bar{w}))|_{z=w}. 
$$
This clearly defines a Wick type deformation quantization of polynomials on $\C^n$, which can be extended to an associative product on formal power series and giving rise to the Wick algebra:
\begin{defn}\label{definition: Wick-algebra}
The {\em Wick product} on the space $\mathcal{W}_{\mathbb{C}^n}:=\mathbb{C}[[z^1,\bar{z}^1,\cdots,z^n,\bar{z}^n]][[\hbar]]$ is defined by 
 \begin{equation*}\label{equation: defn-Wick-product}
  f\star g:=\exp\left(-\hbar\sum_{i=1}^n\frac{\partial}{\partial z^i}\frac{\partial}{\partial \bar{w}^i}\right)(f(z,\bar{z})g(w,\bar{w}))\Big|_{z=w}
 \end{equation*}
 We define a weight on $\mathcal{W}_{\C^n}$ such that $|z^i|=|\bar{z}^j|=1,\ |\hbar|=2$. We let $(\mathcal{W}_{\C^n})_k$ denote those sums of monomials of weight at least $k$. 
\end{defn}

\begin{lem-defn}\label{definition: Fock-representation-of-Wick-algebra}
 The Bargmann-Fock space $\mathcal{F}_{\C^n} := \C[[z^1,\cdots,z^n]][[\hbar]]$ is a representation of $\mathcal{W}_{\C^n}$: we define an action of a monomial $f = z^{\alpha_1}\cdots z^{\alpha_k}\bar{z}^{\beta_1}\cdots\bar{z}^{\beta_l} \in \mathcal{W}_{\C^n}$ on $s\in\mathcal{F}_{\C^n}$ by 
\begin{equation}\label{equation: Toeplitz-operators-C-n}
f\circledast s:=\hbar^l\frac{\partial}{\partial z^{\beta_1}}\circ\cdots\circ\frac{\partial}{\partial z^{\beta_l}}\circ m_{z^{\alpha_1}\cdots z^{\alpha_k}}(s),
\end{equation}
where $m_{z^{\alpha_1}\cdots z^{\alpha_k}}$ denotes the multiplication by $z^{\alpha_1}\cdots z^{\alpha_k}$. It is known that
$$
f\circledast(g\circledast s)=(f\star g)\circledast s,
$$
so this defines an action of the Wick algebra $\mathcal{W}_{\C^n}$ on $\mathcal{F}_{\C^n}$, known as the {\em Bargmann-Fock representation} (or the {\em Wick normal ordering} in physics literature). 
\end{lem-defn}

\subsection{General K\"ahler manifolds: peak sections and large volume limits}
\

For a K\"ahler manifold $X$ whose prequantum line bundle is not trivial,  the situation is much more complicated since there are no polynomial sections as on $\C^n$. The idea is to localize the holomorphic sections of $L^{\otimes k}$ to a neighborhood of any point $x_0\in X$ as an analogue of polynomial sections as in $\C^n$, so that the computation on $\C^n$ can be generalized to $X$.   

First of all, we recall the notion of $K$-frame and $K$-coordinate in K\"ahler geometry. Suppose  for every $x_0\in X$, there exists a local holomorphic frame $e_{L,x_0}$ of $L$ and holomorphic coordinates centered at $x_0$ such that the hermitian inner product of $L$ is locally given by $\langle e_{L,x_0}, e_{L,x_0}\rangle=e^{-\rho_{x_0}}$, where the Taylor expansion of $\rho_{z_0}$ is of the form 
\begin{equation}\label{equation: Taylor-expansion-hermitian-metric}
 		\rho_{x_0}(z,\bar{z})\sim  \delta_{ij}z^i\bar{z}^j+\sum_{|I|,|J|\geq 2}\frac{1}{|I|!|J|!}\frac{\partial^{|I|+|J|}\rho_{x_0}}{\partial z^I\bar{z}^J}(x_0)z^I\bar{z}^J.
 	\end{equation}
 Then $e_{L,x_0}$ and $(z^1,\cdots, z^n)$ are called $K$-frame and $K$-coordinate of order $\infty$ respectively. Throughout this paper, we will assume the existence of these coordinates and frames. 
\begin{rmk}
We refer to \cite{Lu-Shiffman} for more details on the existence and more details of $K$-coordinates and $K$-frames. For instance, complex normal coordinates are $K$-coordinates of order $3$. 
\end{rmk}

The closest analogue of polynomial sections is the notion of peak sections introduced in \cite{Tian}. Peak sections are, roughly speaking, global holomorphic sections of $L^{\otimes k}$ for $k>>0$ whose $L^2$-norm is concentrated around $x_0$ with a prescribed leading term of Taylor expansion at $x_0$. Explicitly, for every multi-index $p=(p_1,\cdots, p_n)$ and $r>|p|=p_1+\cdots+p_n$, there exist normalized peak sections $S_{k,p,r}$ such that, under a $K$-coordinate system and $K$-frame $e_{L,x_0}$, we have 
$$
S_{k,p,r}(z)=\left(z_1^{p_1}\cdots z_n^{p_n}+O(|z|^{2r})\right)\cdot e_{L,x_0}^k;
$$
in this way they give generalizations of polynomial sections in the flat case $X=\C^n$. These peak sections also help us relate holomorphic sections for different tensor powers of $L$. 

Peak sections are closely related to the concept of a {\em large volume limit}. First note that $L^{\otimes k}$ can be seen as the prequantum line bundle on $X$ with the rescaled K\"ahler form $k\omega$. Thus letting the tensor power $k\rightarrow\infty$ is equivalent to taking a large volume limit of $X$. In particular, the unit ball around any point $x_{0}\in ( X,k\omega) $ becomes closer and closer to that in the flat case $0\in \mathbb{C}^{n}$.

To compute the inner product $\langle S_{k,p,r},S_{k,p',r'}\rangle$, since the norms of these peak sections are concentrated in the unit ball $D$ centered at $x_0$, we can focus on the following integral:
$$
\int_D \left(z_1^{p_1}\cdots z_n^{p_n}+O(|z|^{2r})\right)\cdot\overline{\left(z_1^{p'_1}\cdots z_n^{p'_n}+O(|z|^{2r'})\right)}e^{-\frac{\rho_{x_0}(z,\bar{z})}{\hbar}}d\text{vol},
$$
where $\hbar=1/k$. In the large volume limit, the inner products of the peak sections can be approximated by a Gaussian type integral concentrated around $x_0$, and we can use the technique of Feynman graph expansion to compute the asymptotic of the inner products. An obvious difference from the flat case in Section \ref{subsection: flat-space} is that here the volume form of the Gaussian integral contains ``interaction terms'' consisting of the higher order Taylor expansion of $\rho_{x_0}$.

\begin{thm}[Feynman-Laplace]\label{theorem: Feynman-Laplace-Gaussian-integral}
	Let $X$ be a compact $n$-dimensional manifold (possibly with boundary), and let $f$ be a smooth function attaining a unique minimum on $X$ at an interior point $x_0 \in X$, and assume that the Hessian of $f$ is non-degenerate at $x_0$; also, let $\mu=\alpha(x)\cdot e^{g(x)}d^nx$ be a top-degree form. Then the integral 
	$$
	I(\hbar):=\frac{1}{\hbar^{n/2}}\int_X \mu e^{-\frac{1}{\hbar}f(x)}=\frac{1}{\hbar^{n/2}}\int_X \alpha(x)\cdot e^{\frac{-f(x)+\hbar g(x)}{\hbar}}dx_1\cdots dx_n,
	$$
	has the following asymptotic expansion as $\hbar\rightarrow 0^{+}$:
	$$
	I(\hbar)\sim\sum_{k\geq 0}a_k\cdot\hbar^k,
	$$
	where each coefficient $a_k$ is a sum of Feynman weights which depends only on the infinite jets of the functions $f,g$ at the point $x_0$. 
\end{thm}

\subsection{The Hilbert space and Toeplitz action}
\

The idea is to construct the desired Hilbert space by taking the span of these peak sections $\{S_{m,p,r}\}$, similar to the flat case. There is a technical issue: there is the index $r$ in a peak section which measures the order of the error term. Thus we need to consider the following double sequence 
$$
V := \prod_{r\geq 0}\left(\prod_{k\geq 0} H^0(X,L^{\otimes k})\right).
$$
We define a subspace of $V$ consisting of those holomorphic sections which satisfy the following asymptotic properties:
\begin{defn}\label{definition: asymptotic_sequence}
For every point $x_0\in X$, we fix a set $\{S_{k,p,r}\}$ of normalized peak sections centered at $x_0$. A sequence of holomorphic sections $\alpha = \{\alpha_{k,r}\in H^0(X,L^{\otimes k})\}$, regarded as an element in $V$, is called an {\em admissible sequence at $x_0$} if it satisfies the following two conditions:
\begin{enumerate}
 	\item For every fixed $r$, the norm of the sequence $\{\alpha_{k,r}\}_{k>0}$ has a uniform bound: $$\|\alpha_{k,r}\|_k\leq C_r.$$
 	\item There is a sequence of complex numbers $\{a_{p,m}\}_{p,m\geq 0}$ such that, for each fixed $r>0$,  we have
	\begin{equation}\label{equation: asymptotic-general-Kahler}
	\langle\alpha_{k,r}-\sum_{2m+|p|\leq r}a_{p,m}\cdot\frac{1}{k^m}\cdot S_{k,p,r+1},\ S_{k,q,r+1}\rangle_k=O\left(\frac{1}{k^{r+1}}\right),
	\end{equation}
	for any multi-index $q$ with $|q|\leq r$.
\end{enumerate}
We define the subspace $U_{x_0}\subset V$ as the $\mathbb{C}$-linear span of admissible sequences at $x_0$.
\end{defn}
The complex numbers $\{a_{p,k}\}$ are called the {\em coefficients} of the admissible sequence $\alpha$. Note that they are independent of either the tensor power $m$ and the {\em weight index} $r$. The coefficients define a natural equivalence relation $\sim$ on $U_{x_0}$, namely, $\alpha$ is equivalent to $\beta$ (denoted as $\alpha \sim \beta$) if the coefficients of $\alpha-\beta$ are all $0$. The vector space we would like to construct is then simply the quotient by this equivalence relation:
$$
H_{x_0}:=U_{x_0}/\sim.
$$

We have the following explicit isomorphism, by turning the coefficients of an admissible sequence to the coefficients of formal power series. 
\begin{lem}
	We have the following isomorphism of $\mathbb{C}$-vector spaces:
	\begin{equation}\label{equation: asymptotic-sequence-isomorphic-to-Fock-space}
		H_{x_0}\cong \mathbb{C}[[y_1,\cdots, y_n]][[\hbar]].
	\end{equation}
\end{lem}
For every smooth function $f\in C^\infty(X)$, there is a natural action of $f$ on the double sequence $V$ via the corresponding Toeplitz operators: 
\begin{align*}
T_f: \{\alpha_{k,r}\}&\mapsto\{T_{f,k}(\alpha_{k,r})\}.
\end{align*}
\begin{prop}[Lemmas 3.12 and 3.13 in \cite{Chan-Leung-Li}]
For every $f\in C^\infty(X)$, the above action satisfies the following properties:
\begin{enumerate}
 \item Suppose that $\alpha=\{\alpha_{k,r}\}$ is an admissible sequence. Then $\{T_{f,k}(\alpha_{k,r})\}$ is also an admissible sequence for any smooth function $f$.
 \item Suppose that two admissible sequences are equivalent, i.e., $\alpha\sim\beta$. Then for any smooth function $f$, we have $T_f(\alpha)\sim T_f(\beta)$.
\end{enumerate}
\end{prop}
Thus the operator $T_f$ on $V$ preserves the sub-quotient $H_{x_0}$. Moreover, these operators defines a representation of $(C^\infty(X)[[\hbar]],\star_{BT})$ on $H_{x_0}$: 
\begin{thm}[Theorem 3.16 in \cite{Chan-Leung-Li}]\label{theorem: representation-Berezin-Toeplitz-on-asymptotic-sequences}
	Let $x_0\in X$ be any point. The action of $C^\infty(X)[[\hbar]]$ on the vector space $H_{x_0}$ satisfies the following relation:
	\begin{equation}\label{equation: representation-Toeplitz-operator}
		T_f\circ T_g=T_{fg}+\sum_{k\geq 1}\hbar^k\cdot T_{C_k(f,g)}, \hspace{3mm}f,g\in C^\infty(X),
	\end{equation}
	where $C_k(-,-)$ are the bi-differential operators which appear in the Berezin-Toeplitz quantization. 
\end{thm}

Now we give an explicit formula of our representation, in terms of {\em formal Toeplitz operators}. First of all, by turning the real number $\hbar$ in the Feynman-Laplace Theorem to a formal variable, we can define {\em formal integral}: 
\begin{defn}\label{definition: formal_integral}
	For $\phi(y,\bar{y})\in(\W_{\mathbb{C}^n})_3$ and $h(y,\bar{y})\in\W_{\mathbb{C}^n}$, we define the following formal integral:
	$$
	\frac{1}{\hbar^n}\int h(y,\bar{y})\cdot e^{\frac{-|y|^2+\phi(y,\bar{y})}{\hbar}}\in\mathbb{C}[[\hbar]]
	$$
	via the Feynman rule in the Feynman-Laplace Theorem.  
\end{defn}
Here  $h(y,\bar{y})\in\W_{\C^n}$ corresponds to the Taylor expansion of $\alpha(x)$ in Theorem \ref{theorem: Feynman-Laplace-Gaussian-integral}. 
Using this formal integral, we can define a {\em Hilbert space in the formal sense}, namely, its inner product takes values in the formal Laurent series $\mathbb{C}((\sqrt{\hbar}))$:
\begin{defn}
	On the $\mathbb{C}((\sqrt{\hbar}))$-vector space $\mathcal{W}_{\mathbb{C}^n}\otimes_{\mathbb{C}[[\hbar]]}\mathbb{C}((\sqrt{\hbar}))$, we define a complex conjugation by extending the complex conjugation on polynomials in $\mathbb{C}^n$: 
	$$
	(\sqrt{\hbar})^ka_{I,J}y^I\bar{y}^J\mapsto(\sqrt{\hbar})^k\bar{a}_{I,J}\bar{y}^I y^J.
	$$
	Fix $\phi(y,\bar{y})\in(\W_{\mathbb{C}^n})_3$. Then for $f, g\in\mathcal{W}_{\mathbb{C}^n}((\sqrt{\hbar}))$, we define their {\em formal inner product} as the following formal integral:
	\begin{equation}\label{defn:formal-inner-product}
		\langle f, g\rangle:=\frac{1}{\hbar^n}\cdot\int f\bar{g}\cdot e^{\frac{-|y|^2+\phi(y,\bar{y})}{\hbar}},
	\end{equation}
	which is in turn defined using Feynman graph expansions as in Definition \ref{definition: formal_integral} and takes value in $\mathbb{C}((\sqrt{\hbar}))$. 
\end{defn}




With respect to the inner product on the above Hilbert space, there is an orthogonal projection from $\W_{\C^n}$ to the subspace $\C[[z^1,\cdots,z^n]][[\hbar]]$. In particular, for any $f\in \W_{\C^n}$, there is an associated {\em formal Toeplitz operator} defined as the composition of orthogonal projection with multiplication. The following theorem describes these operators in terms of asymptotic of Gaussian integrals:

\begin{thm}[Theorem 2.26 in \cite{Chan-Leung-Li}]\label{theorem: asymptotics-local}
	 Suppose $\varphi(z,\bar{z})$ is a smooth function on $\mathbb{D}^{2n}$ which attains its unique minimum at the origin.  
	Let $f,\varphi,s$ be functions on $\mathbb{D}^{2n}$ such that $\bar{\partial}s=0$, $\varphi$ has a unique minimum at the origin and satisfies 
	\begin{equation}\label{equation: Taylor-expansion-varphi}
	J_{\varphi}=|y|^2+\sum_{I,J\geq 2}\frac{1}{I!J!}\frac{\partial^{|I|+|J|}\varphi}{\partial z^I\partial \bar{z}^J}(0)y^I\bar{y}^J. 
	\end{equation}
	There exist complex numbers $a_{k,I}$ so that for every fixed multi-index $J$, we have the following asymptotics as $\hbar\rightarrow 0$:
	\begin{equation}\label{equation: formula-local-Toeplitz}
		\frac{1}{\hbar^n}\int_{\mathbb{D}^{2n}}\bigg(f\cdot s-\sum_{2k+|I|\leq r} \frac{1}{\hbar^k}a_{k,I}z^I\bigg)\cdot\bar{z}^Je^{-\frac{\varphi(z,\bar{z})}{\hbar}}dvol_{\mathbb{D}^{2n}}=O(\hbar^{r+1}). 
	\end{equation}
	In particular, these $a_{k,I}$'s only depend on the Taylor expansions of $f,s,\varphi$ and $\psi$ at the origin. 
\end{thm}

It follows from asymptotics of inner products of peak sections $S_{m,p,r}$'s as $m\rightarrow \infty$ and equation \eqref{equation: asymptotic-general-Kahler} that $H_{x_0}$ is a formal Hilbert space.  We define $J_{f,x_0}\in\W_{\mathbb{C}^n}$ via the jets of $f$ under $K$-coordinates:
\begin{equation}\label{equation: map-function-to-Wick}
J_{f,x_0}:=\sum_{|I|, |J|\geq 0}\frac{1}{I!J!}\frac{\partial^{|I|+|J|}f}{\partial z^I\bar{z}^J}(x_0)y^I\bar{y}^J,
\end{equation}
where the sum is over all multi-indices.
\begin{thm}[Theorem 3.18 in \cite{Chan-Leung-Li}]\label{proposition: explicit-formula-asymptotic-representation}
	Let $f$ be any smooth function on $X$, and $J_{f,x_0}$ be defined as above. We define $O_{f,x_0}\in\W_{\mathbb{C}^n}$ as the unique solution of the following equation:
	\begin{equation}\label{equation: formal-Toeplitz-I}
	J_{f,x_0}\cdot e^{\Phi/\hbar}=e^{\Phi/\hbar}\star O_{f,x_0}.
	\end{equation}
	Then the action of $T_f$ on $\alpha\in H_{x_0}$ is given by 
	$$
	T_f(\alpha)=O_{f,x_0}\star\alpha.
	$$
	In particular, this implies that the representation $H_{x_0}$ is local in $f\in C^\infty(X)$, i.e., it only depends on the infinite jets of $f$ at $x_0$. 
\end{thm}

\section{Quantization of Kapranov's $L_\infty$ structures on K\"ahler manifolds and Bargmann-Fock sheaves}\label{section: L_infty-structures-and-BF-sheaves}
In this section, we will explain how the the formal Hilbert spaces $H_{x_0}$ constructed in the previous section can be glued together consistently on $X$. A short answer to this question is that there is a sheaf of modules over deformation quantization algebra which admits a flat Fedosov connection. In particular, for every $x_0\in X$, the germ of flat sections at $x_0$ contains a dense subspace of $H_{x_0}$. We will first recall Fedosov connections on K\"ahler manifolds via quantization of $L_\infty$ structure, and explain how these connections can be extended to a module sheaf compatibly.

\subsection{Kapranov's $L_\infty$ structures on K\"ahler manifolds and their quantizations}\label{section: L_infty-structures}
\

Fedosov's connection was introduced in \cite{Fed} to give a geometric construction of deformation quantization on symplectic manifolds. On a K\"ahler manifold $X$, we can implement Fedosov's approach. First of all, there are the following Weyl bundles on $X$:
\begin{align*}\label{equation: Weyl-bundle}
 \W_{X}& := \widehat{\Sym}T^*X[[\hbar]], \quad \overline{\W}_X:=\widehat{\Sym}\overline{T^*X}[[\hbar]],\\
 \W_{X,\C}& := \W_{X}\otimes\overline{\W}_X=\widehat{\Sym}T^*X_{\C}[[\hbar]].
\end{align*}
The fiberwise Hermitian structure on the complexified tangent bundle $TX_{\C}$ enables us to define a fiberwise (non-commutative) Wick product on $\W_{X,\C}$:
\begin{equation*}\label{equation: fiberwise-Wick-product}
 \alpha\star\beta:=\sum_{k\geq 0}\frac{1}{k!}\cdot\left(\frac{\sqrt{-1}\cdot\hbar}{2}\right)^k\omega^{i_1\bar{j}_1}\cdots\omega^{i_k\bar{j}_k}\frac{\partial^k\alpha}{\partial y^{i_1}\cdots\partial y^{i_k}}\frac{\partial^k\beta}{\partial \bar{y}^{j_1}\cdots\partial \bar{y}^{j_k}}.
\end{equation*}
There is the following symbol map:
\begin{equation}\label{equation: symbol-map}
 \sigma: \A_X^\bullet(\W_{X,\C})\rightarrow\A_X^\bullet.
\end{equation}
\begin{defn}
A connection on $\W_{X,\C}$ of the form 
$$
D=\nabla-\delta+\frac{1}{\hbar}[I,-]_{\star}
$$
is called a {\em Fedosov abelian connection} if $D^2=0$.  Here $\nabla$ is the Levi-Civita connection, and $I\in\A^1(X,\W_{X,\C})$ is a $1$-form valued section of $\W_{X,\C}$, with 
\begin{equation}\label{equation: delta}
\delta:=dz^i\wedge\frac{\partial a}{\partial y^i}+d\bar{z}^j\wedge\frac{\partial a}{\partial\bar{y}^j}.
\end{equation}
\end{defn}
Kapranov's $L_\infty$ structure on a K\"ahler manifold can be summarized in the following:
\begin{thm}[Theorem 2.6 and Reformulation 2.8.1 in Kapranov \cite{Kapranov}]\label{thm:Kapranov-L-infinity}
	Let $X$ be a K\"ahler manifold. Then there exist
	\begin{equation}\label{equation: R-n}
	R_n^*\in\mathcal{A}^{0,1}_X(\Hom(T^*X,\Sym^n(T^*X))),\qquad n\geq 2
	\end{equation}
	such that their extensions $\tilde{R}^*_n$ to the holomorphic Weyl bundle $\W_X$ by derivation satisfy
	$$\left(\bar{\partial}+\sum_{n\geq 2}\tilde{R}_n^*\right)^2 = 0,$$
	or equivalently,
\begin{equation}\label{equation: square-zero-0-1-part}
	\bar{\partial}\tilde{R}_n^*+\sum_{j+k=n+1}\tilde{R}_j^*\circ\tilde{R}_k^* = 0
\end{equation}
	for any $n \geq 2$.
\end{thm}
It is easy to show that the $L_\infty$ structure induces a flat connection $D_K$ on the holomorphic component $\W_X$ of the Weyl bundle:
\begin{equation*}\label{eqn:classical-flat-connection-holomorphic}
	D_K=\nabla-\delta^{1,0}+\sum_{n\geq 2}\tilde{R}_n^*. 
\end{equation*}
Here $\delta^{1,0}$ denotes the $(1,0)$-part of $\delta$ in equation \eqref{equation: delta}. The symbol map induces a natural one-to-one correspondence $\sigma:\Gamma^{flat}(U,\W_X)\cong \mathcal{O}_X(U)$. For a (local) holomorphic function $f$, the associated flat section under $D_K$ gives the Taylor expansion of $f$ under $K$-coordinate at every point $x_0\in X$. And we denote this flat section by $J_f$ (``J'' for jets).

\subsubsection{Classical and quantum extensions of $D_K$}
\

There are two ways to extend $D_K$ to flat connections on the complexified Weyl bundle $\W_{X,\C}$, one classical and the other quantum.  For the classical extension, notice that the complex conjugate of the connection $D_K$ is a flat connection $\overline{D}_K$ on $\overline{\mathcal{W}}_X$. Then
$$D_C := D_K \otimes 1+1\otimes\overline{D}_K$$
is naturally a flat connection on $\mathcal{W}_{X,\mathbb{C}}=\W_{X}\otimes\overline{\W}_X$ such that $D_C|_{\W_X}=D_K$. The geometry behind this extension is very simple: since the flat sections with respect to $D_K$ correspond to (local) holomorphic functions on $X$, by adding the anti-holomorphic components in $\mathcal{W}_{X,\mathbb{C}}$, we shall see all the smooth functions. 
This is indeed the case.
\begin{prop}[Proposition 2.10 in \cite{CLL3}]\label{proposition: classical-flat-section-smooth-function}
There is a one-to-one correspondence between $C^\infty(X)[[\hbar]]$ and the space of flat sections of the Weyl bundle $\W_{X,\mathbb{C}}$ with respect to the flat connection $D_C$. 
\end{prop}
We also denote by the flat section associated to a smooth function by $J_f$ since it gives the Taylor expansion of $f$ under $K$-coordinates and their conjugates. This is also the reason we call $D_C$ a classical extension of $D_K$. 

In \cite{CLL}, we constructed a family of Fedosov abelian connections on K\"ahler manifolds via a quantum extension of $D_K$. Explicitly, we use the following $\A_X^\bullet$-linear operator 
$$
L:\A_X^\bullet\left(\widehat{\Sym}(T^*X)\otimes TX\right)\rightarrow \A_X^\bullet\left(\widehat{\Sym}(T^*X)\otimes \overline{T^*X}\right)
$$
``lifting the last subscript'' using the K\"aher form to define 
$$
I_n:=L(R_n^*)=R_{i_1\cdots i_n,\bar{l}}^j\omega_{j\bar{k}}d\bar{z}^l\otimes (y^{i_1}\cdots y^{i_n}\bar{y}^{k})\in\mathcal{A}_X^{0,1}(\mathcal{W}_{X,\mathbb{C}}). 
$$
There is a simple observation behind this operator $L$: we simply replace the contraction between $TX$ in $\tilde{R}_n$ and $T^*X$ in $\W_X$ by a bracket $[-,-]_\star$, which can be extended to $\W_{X,\C}$. 
Let $I:=\sum_{n\geq 2}I_n$. Then we proved in \cite{CLL} that 
\begin{thm}[Theorems 2.17 and 2.25 in \cite{CLL}]
Let $\alpha$ be a representative of a formal cohomology class in $\hbar H^2_{dR}(X)[[\hbar]]$ of type $(1,1)$. Then there exists a solution of the Fedosov equation of the form $I_\alpha = I+ J_\alpha \in \mathcal{A}_X^{0,1}(\mathcal{W}_{X,\mathbb{C}})$:
 \begin{equation}\label{equation: Fedosov-equation}
\nabla I_\alpha - \delta I_\alpha + \frac{1}{\hbar} I_\alpha\star I_\alpha + R_\nabla=\alpha.
\end{equation}
We denote the corresponding Fedosov abelian connection by $D_{F,\alpha}$. Every such $D_{F,\alpha}$ is a quantum extension of $D_K$, i.e., $D_{F,\alpha}|_{\W_X}=D_K$. And the deformation quantization associated to $D_{F,\alpha}$ is a Wick type star product whose Karabegov form is $2\sqrt{-1}\omega-\alpha$. 
\end{thm}
There are  several nice properties of these Fedosov connections $D_{F,\alpha}$. First of all, 
the Karabegov form of the associated star products is clear from the connection. Secondly, the connection $D_{F,\alpha}$ ``looks classical'' because it does {\em not} contain any $\hbar$. Lastly, the equality $D_{F,\alpha}|_{\W_X}=D_K$ implies that for a (local) holomorphic function $f$, we have $D_{F,\alpha}(J_f)=0$. This is saying that holomorphic functions do not receive any ``quantum corrections''. This property makes the Fedosov quantization more compatible with the Berezin-Toeplitz quantization: when $f$ is holomorphic, the Toeplitz operator $T_{f,k}=\Pi_k\circ m_f$ is simply the ``classical'' operator given by multiplication by $f$.

\subsection{Bargmann-Fock sheaves}\label{section: BF-sheaves}
\

The idea behind Fedosov's approach to deformation quantization is clear: the fiberwise Moyal-Weyl product on the Weyl bundle describes the local picture of star products, and a Fedosov abelian connection gives the gluing data for a global star product. We follow the same line of thought in our construction of Bargmann-Fock sheaves on K\"ahler manifolds, starting with the fiberwise Bargmann-Fock action of the Weyl bundle $\W_{X,\C}$ (equipped with the fiberwise Wick product $\star$) on $\W_X=\widehat{\Sym} T^*X$ as in Definition \ref{definition: Fock-representation-of-Wick-algebra}. 

It is easy to see that there is a unique choice of a connection on $\W_X[\hbar^{-1}]$ compatible with $D_{F,\alpha}$ in the obvious sense, which is explicitly given by 
\begin{equation}\label{equation: connection-D-alpha}
D_\alpha=\nabla+\frac{1}{\hbar}\gamma_\alpha\circledast-.
\end{equation}
However, this connection is not flat:
\begin{lem}[Lemma 3.8 in \cite{CLL3}]\label{lemma: curvature-Levi-Civita-extended-holomorphic-Weyl}
	The curvature of $D_\alpha$ is given by
	$$D_\alpha^2 = \frac{1}{\hbar}\omega_\hbar - \text{Ric}_X,$$
	where $\omega_\hbar=2\sqrt{-1}\omega-\alpha$ and $\text{Ric}_X = R_{i\bar{j}k}^kdz^i\wedge d\bar{z}^j$ is the Ricci form of $X$.
	In particular, the connection $D_\alpha$ on $\W_{X,e}$ is not flat.
\end{lem}
A naive idea is to twist $\W_X$ by a line bundle to cancel its curvature. However, the curvature $(1,1)$-form $\frac{1}{\hbar}\omega_\hbar$ involves the formal variable $\hbar$ and does not satisfy the integrality condition. What we need is the notion of {\em formal line bundles}, which we formulate as sheaves. First of all, we consider the following extensions of $\W_X$ and $\mathcal{O}_X$ by allowing formal exponentials:
\begin{defn}\label{definition: Weyl-bundle-with-exponentials}
We define the sheaf $\W_{X,e}$ of extended  Weyl algebra with exponentials as follows:
for every open set $U\subset X$, we consider the space of finite sum of pairs
\begin{equation}\label{equation: Weyl-bundle-with-exponentials}
\sum_{i=1}^k(f_i, e^{g_i/\hbar}),
\end{equation}
where $f_i,g_i$'s are smooth sections of $\W_X$ on $U$. We define the multiplication by the linear extension of 
$$
(f_1,e^{g_1/\hbar})\cdot(f_2,e^{g_2/\hbar}):=(f_1f_2,e^{(g_1+g_2)/\hbar}).
$$
These are subject to the equivalence relation that $(f_1,e^{g_1/\hbar})\sim (f_2,e^{g_2/\hbar})$ if $f_1=f_2$ and $g_1-g_2\in\C[[\hbar]]$. Then the space $\W_{X,e}(U)$ of sections of $\W_{X,e}$ over $U$ is given by the set of equivalence classes. There is a sub-sheaf $\mathcal{O}_{X,e}$ of $\W_{X,e}$ consisting of those finite sums in equation \eqref{equation: Weyl-bundle-with-exponentials} where $f_i,g_i\in\mathcal{O}_X(U)[[\hbar]]$
\end{defn}
\begin{notn}
We will use the notation $f\cdot e^{g/\hbar}$ for the the pair $(f,e^{g/\hbar})$ in $\W_{X,e}$ (and also $\mathcal{O}_{X,e}$). 
\end{notn}
The fiberwise Bargmann-Fock action can be extended on $\W_{X,e}$ with some care: for a monomial as in Definition \ref{definition: Fock-representation-of-Wick-algebra}, we can extend its action to $\W_{X,e}$ by the same differential operator as in equation \eqref{equation: Toeplitz-operators-C-n}. In particular, 
$$
\bar{y}^j\circledast(f\cdot e^{g/\hbar})=\hbar\frac{\omega^{i\bar{j}}}{2\sqrt{-1}}\frac{\partial}{\partial y^i}(f\cdot e^{g/\hbar})=\hbar\frac{\omega^{i\bar{j}}}{2\sqrt{-1}}\left(\frac{\partial f}{\partial y^i}+\frac{1}{\hbar}f\cdot\frac{\partial g}{\partial y^i}\right)\cdot e^{g/\hbar}
$$
However, for general elements of $\W_{X,\C}$ on $\W_{X,e}$ we could run into infinite sums such as the following example: when $X=\C$, $g=y$ and $f=\sum_{k\geq 1}\bar{y}^k$, 
\begin{align*}
f\circledast e^{g/\hbar}=\left(\sum_{k\geq 1}\bar{y}^k\right)\circledast e^{y/\hbar}
=\sum_{k\geq 1}\left(\hbar\partial_{y}\right)^k (e^{y/\hbar})
=e^{y/\hbar}+e^{y/\hbar}+\cdots.
\end{align*}
If we write $f=\sum_{k,I,J}\hbar^k f_{I,\bar{J},k}y^I\bar{y}^J$ and $g=\sum_{I,k}\hbar^k g_{k,I}y^I$, then it is not difficult to see from the above example that the infinite sums come from two sources:
\begin{enumerate}
 \item Those terms $g_{0,i}y^i$ in $g$ which are linear in $\W_{X}$ and do not include $\hbar$;
 \item The infinite sums $\sum_{J}f_{k_0,I_0,\bar{J}}\hbar^{k_0}y^{I_0}\bar{y}^J$ for fixed indices $I_0, k_0$. 
\end{enumerate}
A section $\alpha = \sum_{k,I,J}\hbar^k\alpha_{k,I,\bar{J}}y^I\bar{y}^J$ of $\W_{X,\C}$ is called {\em admissible} if it satisfies the condition that for every fixed $I_0$ and $k_0$, $\sum_{J}\hbar^{k_0}\alpha_{k_0,I_0,\bar{J}}y^{I_0}\bar{y}^J$ is a finite sum.
It is easy to see from the construction that the term $\gamma_\alpha$ in equation \eqref{equation: connection-D-alpha} is admissible. (We refer to the proof of Theorem 2.17 in \cite{CLL} for the details). Thus $D_\alpha$ extends to a well-defined connection on $\W_{X,e}$ whose curvature is $\frac{1}{\hbar}\omega_\hbar - \text{Ric}_X$.

We now define a formal line bundle as an invertible $\mathcal{O}_{X,e}$-module, on which we can also define connection and curvature. This is analogue to the definition of holomorphic line bundle: the exponentials in $\mathcal{O}_{X,e}$ plays the role of the transition function. They are also similar to local line bundles in \cite{Melrose} and twisting bundles in \cite{Tsygan}.
\begin{prop}[Lemma 3.5 in \cite{CLL3}]\label{proposition: formal-holomorphic-line-bundle-correspondence-1-1-class}
For every formal closed $(1,1)$-form $\alpha\in \A_{closed}^{1,1}(X)[[\hbar]]$, there is a formal line bundle $L_{\alpha/\hbar}$ with connection $\nabla_{L_{\alpha/\hbar}}$ whose curvature is given by $\frac{1}{\hbar}\cdot\alpha$.  
\end{prop}
Thus we can take the tensor product of $\W_{X,e}$ with a  formal line bundle, on which there exists a compatible Fedosov flat connection.
\begin{defn}\label{definition: Bargmann-Fock-sheaf}
	For a representative $\alpha$ of a formal $(1,1)$-class $[\alpha]\in \hbar H^{1,1}_{dR}(X)[[\hbar]]$, let $\omega_\hbar := 2\sqrt{-1}\cdot\omega-\alpha$ and $\alpha' := -\omega_\hbar+\hbar\cdot\text{Ric}_X$. Then we define the {\em sheaf of Bargmann-Fock modules} as $$\mathcal{F}_{X,\alpha}:=\W_{X,e}\otimes_{\mathcal{O}_{X,e}} L_{\alpha'/\hbar}.$$
	It is equipped with the connection
	$$D_{B,\alpha} := (\nabla+\frac{1}{\hbar}\gamma_{\alpha}\circledast -)\otimes \nabla_{L_{\alpha'/\hbar}}.$$
\end{defn} 

\begin{lem}[Lemma 3.10 in \cite{CLL3}]\label{lemma: compatibility-Fedosov-connection} 
The connection $D_{B,\alpha}$ is compatible with the Fedosov connection $D_{F,\alpha}$, i.e., if $O\in\W_{X,\C}$ is an admissible section and $s$ is a section of $\mathcal{F}_{X,\alpha}$ , then we have
\begin{equation*}\label{equation: compatibility-Fedosov-connection}
D_{B,\alpha}(O\circledast s)=D_{F,\alpha}(O)\circledast s+(-1)^{|O|} O\circledast(D_{B,\alpha}(s)).
\end{equation*}
In particular, if $O$ and $s$ are flat sections of $\W_{X,\C}$ and $\mathcal{F}_{X,\alpha}$ respectively, then $O \circledast s$ is also a flat section.
\end{lem}

Since formal smooth functions $C^\infty(X)[[\hbar]]$ is canonically isomorphic to the flat sections of $\W_{X,\C}$ under $D_{F,\alpha}$, it is natural to expect that the above compatibility of $D_{B,\alpha}$ and $D_{F,\alpha}$ implies that those flat sections of $\mathcal{F}_{X,\alpha}$ form a module over $C_X^\infty(X)[[\hbar]]$.  However, we do not expect that the flat section $O_f$ corresponding to a smooth function to be admissible. Thus an infinite sum can not be avoided, and we need to assume certain analytic condition on both functions and sections of $\mathcal{F}_{X,\alpha}$. 

We define the following sub-class of real analytic functions, which roughly speaking consists of analytic functions with analyticity at least the same as that of $\omega_\hbar$. Precisely, for any $x_0\in X$ we let $r(x_0)$ denote the radius of convergence of $\omega_\hbar$ under a $K$-coordinate centered at $x_0$.
\begin{lem-defn}\label{definition: convergence-property-analytic-function}
For every open set $U\subset X$, let $C^{\omega_\hbar}_X(U)[[\hbar]]$ denote the set of real analytic functions on $X$ such that at every point $x_0 \in X$, the radius of convergence is greater than or equal to $r(x_0)$ under a $K$-coordinate centered at $x_0$. These functions are closed under the star product $\star_\alpha$, and form a sheaf $C^{\omega_\hbar}_X[[\hbar]]$ of algebras on $X$ under the Fedosov star product $\star_\alpha$. 
\end{lem-defn}

\begin{defn}\label{definition: flat-Bargmann-Fock}
The {\em Bargmann-Fock sheaf} $\mathcal{F}_{X,\alpha}^{\text{flat}}$ is defined as the sub-sheaf of $\mathcal{F}_{X,\alpha}$ which consists of flat sections that are finite sums of the following form: $\alpha\cdot e^{\beta/\hbar}\otimes e_U$, where we can write $\beta=\sum_{|I|\geq 0}\beta_Iy^I$ locally.
We require that the coefficients of the degree $1$ terms, i.e., $\beta_i$, $1\leq i\leq n$, satisfy the following boundedness condition:
\begin{equation}\label{equation: convergent-Fock-boundedness-condition}
\left\|\sum_{i=1}^n\beta_iy^i\right\|_{x_0} < r(x_0),
\end{equation}
where the norm is defined using the Hermitian metric on $T^*X$. 
\end{defn}


These analytic conditions guarantee that for $f\in C^{\omega_\hbar}_X(U)$ and $s\in\mathcal{F}_{X,\alpha}^{\text{flat}}(U)$, there is a well-defined $O_f\circledast s\in\mathcal{F}_{X,\alpha}^{\text{flat}}(U)$. 
\begin{thm}[Theorem 3.17 in \cite{CLL3}]\label{theorem: Bargmann-Fock-representation}
 The Bargmann-Fock sheaf $\mathcal{F}_{X,\alpha}^{\text{flat}}$ is a sheaf of modules over $\left(C^{\omega_\hbar}_X[[\hbar]],\star_\alpha\right)$. 
\end{thm}

Our construction and proof of Theorem \ref{theorem: Bargmann-Fock-representation} in \cite{CLL3} are analytic in nature and follow Fedosov's original approach closely.
On the other hand, closely related studies on such module sheaves have been carried out algebraically using deformation-obstruction theory. In the real symplectic manifolds context, such constructions were established in the work of Nest-Tsygan \cite{Nest-Tsygan} and Tsygan \cite{Tsygan}. In \cite{BGKP}, Baranovsky, Ginzburg, Kaledin and Pecharich gave a deformation theoretic construction of quantizations of line bundles as module sheaves in the algebraic setting. 

 \subsubsection{Prequantum Bargmann-Fock sheaf and Berezin-Toeplitz quantization}
 \

 In particular, when the K\"ahler manifold $X$ is prequantizable and $\alpha = \hbar c_1(X)$, it is known that the star product is exactly given by the Berezin-Toeplitz star product $\star_{BT}$ and the corresponding formal line bundle
precisely characterizes the asymptotics of the tensor powers $L^{\otimes k}$ of the prequantum line bundle as $k\rightarrow\infty$. 
 
In this subsection, we consider the example of $\mathcal{F}_{X,\alpha}$ where $\alpha = -\hbar \text{Ric}_X =  -\hbar\cdot R_{i\bar{j}k}^kdz^i\wedge d\bar{z}^j$. Explicitly,
$$
\mathcal{F}_{X, \alpha} = \W_{X,e}\otimes_{\mathcal{O}_{X,e}}L_{-2\sqrt{-1}\omega/\hbar},
$$
We choose a holomorphic frame $e_{x_0}$ of $L_{-2\sqrt{-1}\omega/\hbar}$ satisfying the following condition:
$$
\nabla_{L_{-2\sqrt{-1}\omega/\hbar}}(e_{x_0})=-\frac{1}{\hbar}\partial\rho_{x_0}\otimes e_{x_0},
$$
and define a local section of  $\W_X$ by $\beta=\sum_{k\geq 1}(\tilde{\nabla}^{1,0})^k(\rho_{x_0})$.
\begin{thm}[Theorem 4.4 in \cite{CLL3}]\label{theorem: flat-section-Bargmann-Fock}
There is $e^{\beta/\hbar}\otimes e_{x_0}\in\mathcal{F}_{X,\alpha}^{\text{flat}}$. And a section of $\mathcal{F}_{X,\alpha}$ of the form $A\cdot e^{\beta/\hbar}\otimes e_{x_0}$ around $x_0$ is flat, i.e., $D_{B,\alpha}\left(A\cdot e^{\beta/\hbar}\otimes e_{x_0}\right)=0$ if and only if $D_K(A)=0$, or equivalently, $A=J_s$ for some holomorphic function $s$. 
\end{thm}

According to this theorem, there exists a subspace $V_{x_0}$ of the stalk $(\mathcal{F}_{X,\alpha}^{\text{flat}})_{x_0}$, which is isomorphic to the space of germs of formal holomorphic functions at $x_0$ , i.e., $V_{x_0}\cong\mathcal{O}_{X,x_0}[[\hbar]]$, such that  $V_{x_0}$ is a representation of the Berezin-Toeplitz deformation quantization. 

After a choice of a $K$-coordinate centered at $x_0$, the vector space $V_{x_0}$ is naturally a subspace of $H_{x_0}$ consisting of those formal power series which are convergent in some neighborhood of $0\in\C^n$. 
\begin{thm}[Theorem 4.10 in \cite{CLL3}]\label{theorem: Bargmann-Fock-microlocal-Toeplitz}
The representation of $C^{\omega_\hbar}(X)[[\hbar]]$ on $V_{x_0}$ is given explicitly as follows: 
Let $f\in C^{\omega_\hbar}(X)[[\hbar]] $ and $\Psi_s:=J_s\cdot e^{\beta/\hbar}\otimes e_{x_0}\in V_{x_0}\subset(\mathcal{F}_{X,\alpha}^{\text{flat}})_{x_0}$ where $s\in\mathcal{O}_{X,x_0}[[\hbar]]$. Then 
\begin{equation}\label{equation: formal-Toeplitz}
O_f\circledast\left(J_s\cdot e^{\beta/\hbar}\otimes e_{x_0}\right)=J_{s'}\cdot e^{\beta/\hbar}\otimes e_{x_0},
\end{equation}
where $s'\in\mathcal{O}_{X,x_0}[[\hbar]]$ is determined by its jets $J_{s'}$ at $x_0$ explicitly given by
$$
T_{(J_f)_{x_0},\Phi}\left(J_s\right) = J_{s'}.
$$
\end{thm}
By comparing equation \eqref{equation: formal-Toeplitz-I} with \eqref{equation: formal-Toeplitz}, we can see that  $s^{\prime }$ is obtained by applying the formal Toeplitz operator associated to $J_f$ to $J_s$. Thus this theorem tells us that, modulo technical issues (dense subspace, analyticity), every representation $H_{x_0}$ is a subspace of the stalk $(\mathcal{F}_{X,\alpha})_{x_0}$ of Bargmann-Fock sheaf. 

\section{BV quantization of K\"ahler manifolds}\label{section: BV}

In this section we give an application of the Fedosov connection. In \cite{GLL}, it is shown that the physical interpretation of deformation quantization and algebraic index theorem of a symplectic manifold $X$ is the Batalin-Vilkovisky (BV) quantization of the one-dimensional Chern-Simons theory, which describes a sigma model with domain $S^1$ and target $X$ in a neighborhood of the constant maps. This model is quantized rigorously using Costello's formalism of effective renormalization \cite{Kevin-book}. An effective quantization must satisfy the {\em Quantum Master Equation} (QME) which describes the gauge invariance at the quantum level. In this model, the QME can be described using the geometry of BV bundles:
\begin{defn}[cf. Definition 2.19 in \cite{GLL}]\label{definition: BV-bundle}
 The {\em BV bundle} of a K\"ahler manifold $X$ is defined to be 
 $$
 \widehat{\Omega}^{-\bullet}_{TX}:=\widehat{\Sym}(T^*X_{\C})\otimes\wedge^{-\bullet}(T^*X_{\C}),\quad \wedge^{-\bullet}(T^*X_{\C}):=\bigoplus_k\wedge^k(T^*X_{\C})[k],
 $$
 where $\wedge^k(T^*X_{\C})$ has cohomological degree $-k$.
\end{defn}
\begin{defn}[Definition 2.22 in \cite{GLL}]\label{definition: QME}
The operator 
$$Q_{BV}:=\nabla+\hbar\Delta+\frac{1}{\hbar}d_{TX}R_{\nabla}$$
is a differential on the BV bundle (i.e., $Q_{BV}^2=0$), which we call the {\em BV differential}. 
A section $\gamma_\infty$ of the BV bundle is said to {\em satisfy the quantum master equation (QME)} if
\begin{equation}\label{equation: quantum-master-equation}
 Q_{BV}(e^{\gamma_\infty/\hbar})=0.
\end{equation}
\end{defn}


A solution $\gamma_\infty$ of the QME \eqref{equation: quantum-master-equation} induces a differential 
\begin{equation}\label{equation: differential-global-quantum-observable}
\nabla + \hbar\Delta + \{\gamma_\infty,-\}_\Delta
\end{equation}
on the BV bundle.
The cochain complex $(\A_X^\bullet\left(\hat{\Omega}_{TX}^{-\bullet}\right)[[\hbar]], \nabla  + \hbar\Delta + \{\gamma_\infty,-\}_\Delta)$ is called the {\em global quantum observables} of the one-dimensional Chern-Simons model. In \cite{GLL}, it is shown that the {\em fiberwise Berezin integration}, defined by taking the top degree component in odd variables and setting the even variables to $0$:
\begin{equation*}\label{equation: Berezin-integration}
\int_{Ber}:\A_X^\bullet\left(\hat{\Omega}_{TX}^{-\bullet}\right) \rightarrow \A_X^\bullet,
\qquad a \mapsto\frac{1}{n!}(\iota_{\Pi})^n(a)\bigg|_{y^i=\bar{y}^j=0},
\end{equation*}
is a cochain map, with respect to the BV differential $Q_{BV}$ on $\A_X^\bullet\left(\hat{\Omega}_{TX}^{-\bullet}\right)$ and the de Rham differential on $\A_X^\bullet$.   
Thus we get a well-defined composition map on cohomology classes, which enables us to define the correlation functions (or expectation values) of global quantum observables:
\begin{equation*}\label{equation: BV-integration-composite-ordinary-integration}
H^*(\A_X^\bullet\left(\hat{\Omega}_{TX}^{-\bullet}\right)[[\hbar]])\overset{\int_{Ber}}{\longrightarrow} H^*_{dR}(X)[[\hbar]]\overset{\int_{X}}{\longrightarrow}\mathbb{C}[[\hbar]].
\end{equation*}


The main ingredient in BV quantization is to find a solution for the Quantum Master Equation (QME). This can be obtained by applying the {\em homotopy group flow operator} using the propagator $\partial_P$ of the one-dimensional sigma model to a solution of the Fedosov equation. Our formulation here follows that in \cite{GLL}*{Section 2.4} but with a significant modification of the definition of the propagator in order to adapt to the K\"ahler setting. We refer to \cites{CLL,GLL} for details on the explanation of the notations here. 
\begin{thm}[Theorem 2.26 in \cite{GLL} and Theorem 3.15 in \cite{CLL}]\label{theorem: Fedosov-connection-RG-flow-QME}
Let $\gamma\in\A_X^1(\W_{X,\C})$, we define $\gamma_\infty\in\A_X^\bullet(\hat{\Omega}_{TX}^{-\bullet})[[\hbar]]$ by 
 $$
 e^{\gamma_{\infty}/\hbar}:=\text{Mult}\int_{S^1[*]}e^{\hbar\partial_P+D}e^{\otimes\gamma/\hbar}.
 $$
Suppose that $\gamma$ is a solution of the Fedosov equation. 
Then $e^{\tilde{R}_\nabla/2\hbar}\cdot e^{\gamma_\infty/\hbar}$ is a solution of the QME \eqref{equation: quantum-master-equation}, i.e., $Q_{BV}\left(e^{\tilde{R}_\nabla/2\hbar}\cdot e^{\gamma_\infty/\hbar}\right)=0$.
\end{thm}

It turns out that, in the K\"ahler case, if $\gamma$ is the solution of the Fedosov equation  obtained by quantizing Kapranov's $L_\infty$-algebra structure, then the Feynman graph expansion of the QME solution $\gamma_\infty$ involves {\em only trees and one-loop graphs}; in other words, $\gamma_\infty$ gives a {\em one-loop exact BV quantization} of the K\"ahler manifold $X$. This is in sharp contrast with the general symplectic case \cite{GLL}, in which the BV quantization involves all-loop quantum corrections.
The same kind of one-loop exactness was observed for the holomorphic Chern-Simons theory by Costello \cite{Kevin-CS} and for a sigma model from $S^1$ to the target $T^*Y$ (cotangent bundle of a smooth manifold $Y$) by Gwilliam-Grady \cite{Gwilliam-Grady}.

In the BV formalism, there is a rich structure of {\em factorization algebra} of quantum observables, including the local-to-global factorization map of quantum observables. In our models, the cochain complex of local quantum observable is $\left(\A_X^{\bullet}(\mathcal{W}_{X,\C}),D_{F,\alpha}\right)$.  And the {\em factorization map} is explicitly
\begin{align*}
[-]_\infty: \A_X^{\bullet}(\mathcal{W}_{X,\C})&\rightarrow\A_X^\bullet(\hat{\Omega}_{TX}^{-\bullet})[[\hbar]]\\
O&\mapsto [O]_\infty:=e^{-\gamma_\infty/\hbar}\cdot \left(\text{Mult}\int_{S^1[*]}e^{\hbar\partial_P+D}(O d\theta_1\otimes e^{\otimes\gamma/\hbar})\right).
\end{align*}
This enables us to define the correlation function of a local quantum observable $O_f$ corresponding to a smooth function $f\in C^\infty(X)$ as 
$$
\langle f\rangle:= \int_X\circ\int_{Ber}[O_f]_\infty.
$$
A computation shows that $\langle f \star g\rangle=\langle g \star f\rangle$, and that $\langle f\rangle=\int_X f\cdot\omega^n+O(\hbar)$. Thus the correlation function gives rise to the (unique) {\em trace} of the deformation quantization algebra, as in the following
\begin{defn}
	Let $(C^\infty(X)[[\hbar]],\star)$ denote a deformation quantization of $X$. A {\em trace} of the star product $\star$ is a linear map $\Tr:C^\infty(X)[[\hbar]]\rightarrow\C[[\hbar]]$ such that 
	\begin{enumerate}
		\item $\Tr(f\star g)=\Tr(g\star f)$;
		\item $\Tr(f)=\int_X f\cdot\omega^n+O(\hbar)$.
	\end{enumerate}
	In particular, $\Tr(1)$ is called the {\em algebraic index} of $\star$. 
\end{defn}

As a corollary, we obtain a succinct explicit expression of the algebraic index $\Tr(1)$ of the deformation quantization algebra: 
\begin{thm}[Theorem 3.39 and Corollary 3.40 in \cite{CLL}]\label{theorem:algebraic-index-theorem}
	Let $\gamma$ be a solution of the Fedosov equation and $\gamma_\infty$ be the associated solution of the QME as defined in Theorem \ref{theorem: Fedosov-connection-RG-flow-QME}. Then we have
$$\sigma\left(e^{\hbar\iota_{\Pi}} (e^{\tilde{R}_\nabla/2\hbar}e^{\gamma_\infty/\hbar}) \right)
= \hat{A}(X)\cdot e^{-\frac{\omega_\hbar}{\hbar}+\frac{1}{2}\Tr(\mathcal{R}^+)} = \text{Td}(X)\cdot e^{-\frac{\omega_\hbar}{\hbar}+\Tr(\mathcal{R}^+)},$$
where $Td(X)$ is the Todd class of $X$. In particular, the trace of the function $1$ is given by 
	\begin{align*}
	\Tr(1) = \int_X\hat{A}(X)\cdot e^{-\frac{\omega_\hbar}{\hbar}+\frac{1}{2}\Tr(\mathcal{R}^+)}
	=  \int_X\text{Td}(X)\cdot e^{-\frac{\omega_\hbar}{\hbar}+\Tr(\mathcal{R}^+)}.
	\end{align*}
\end{thm}
This is a cochain level enhancement of the result in \cite{GLL}: there the technique of equivariant localization was applied to show that all graphs of higher genera ($\geq 2$) give rise to exact differential forms after the Berezin integration and thus do not contribute after integration over $X$, while the Feynman weights associated to these graphs in our QME solution $\gamma_\infty$ vanish already on the cochain level.



\begin{bibdiv}
\begin{biblist}

\bib{BGKP}{article}{
    AUTHOR = {Baranovsky, V.},
    author = {Ginzburg, V.},
    author = {Kaledin, D.},
    author = {Pecharich, J.},
     TITLE = {Quantization of line bundles on lagrangian subvarieties},
   JOURNAL = {Selecta Math. (N.S.)},
    VOLUME = {22},
      YEAR = {2016},
    NUMBER = {1},
     PAGES = {1--25},
}

\bib{Bordemann}{article}{
    AUTHOR = {Bordemann, M.},
    author = {Waldmann, S.},
     TITLE = {A {F}edosov star product of the {W}ick type for {K}\"{a}hler
              manifolds},
   JOURNAL = {Lett. Math. Phys.},
    VOLUME = {41},
      YEAR = {1997},
    NUMBER = {3},
     PAGES = {243--253},
}

\bib{Bordemann-Waldmann}{article}{
    AUTHOR = {Bordemann, M.},
    author = {Waldmann, S.},
     TITLE = {Formal {GNS} construction and states in deformation
              quantization},
   JOURNAL = {Comm. Math. Phys.},
    VOLUME = {195},
      YEAR = {1998},
    NUMBER = {3},
     PAGES = {549--583},
}

\bib{Bordemann-Meinrenken}{article}{
    AUTHOR = {Bordemann, M.},
    author = {Meinrenken, E.},
    author = {Schlichenmaier, M.},
     TITLE = {Toeplitz quantization of {K}\"{a}hler manifolds and {${\rm
              gl}(N)$}, {$N\to\infty$} limits},
   JOURNAL = {Comm. Math. Phys.},
    VOLUME = {165},
      YEAR = {1994},
    NUMBER = {2},
     PAGES = {281--296},
}


\bib{Chan-Leung-Li}{article}{
   author={Chan, K.},
   author={Leung, N. C.},
   author={Li, Q.},
   title={A geometric construction of representations of the Berezin-Toeplitz quantization},
   eprint={arXiv:2004.00523 [math-QA]},
}

\bib{CLL}{article}{
   author={Chan, K.},
   author={Leung, N. C.},
   author={Li, Q.},
   title={Kapranov's $L_\infty$ structures, Fedosov's star products, and one-loop exact BV quantizations on K\"ahler manifolds},
   eprint={ arXiv:2008.07057 [math-QA]},
}

\bib{CLL3}{article}{
   author={Chan, K.},
   author={Leung, N. C.},
   author={Li, Q.},
   title={Bargmann-Fock sheaves on K\"ahler manifolds},
   eprint={ arXiv:2008.11496 [math-DG]},
}


\bib{Kevin-book}{book}{
   author={Costello, K.},
   title={Renormalization and effective field theory},
   series={Mathematical Surveys and Monographs},
   volume={170},
   publisher={American Mathematical Society},
   place={Providence, RI},
   date={2011},
   pages={viii+251},
   isbn={978-0-8218-5288-0},
}

\bib{Kevin-CS}{article}{
   author={Costello, K.},
   title={A geometric construction of the Witten genus, II},
   eprint={arXiv:1111.4234 [math.QA]},
  }


\bib{DLS}{article}{
     author= {Dolgushev, V. A. },
     author= {Lyakhovich, S. L.},
     author= {Sharapov, A. A.},
     TITLE = {Wick type deformation quantization of {F}edosov manifolds},
   JOURNAL = {Nuclear Phys. B},
    VOLUME = {606},
      YEAR = {2001},
    NUMBER = {3},
     PAGES = {647--672},
}


\bib{Donaldson}{incollection}{
	AUTHOR = {Donaldson, S. K.},
	TITLE = {Planck's constant in complex and almost-complex geometry},
	BOOKTITLE = {X{III}th {I}nternational {C}ongress on {M}athematical
		{P}hysics ({L}ondon, 2000)},
	PAGES = {63--72},
	PUBLISHER = {Int. Press, Boston, MA},
	YEAR = {2001},
}


\bib{Fed}{article}{
    AUTHOR = {Fedosov, B. V.},
     TITLE = {A simple geometrical construction of deformation quantization},
   JOURNAL = {J. Differential Geom.},
    VOLUME = {40},
      YEAR = {1994},
    NUMBER = {2},
     PAGES = {213--238}
}


\bib{Fedbook}{book}{
    AUTHOR = {Fedosov, B. V.},
     TITLE = {Deformation quantization and index theory},
    SERIES = {Mathematical Topics},
    VOLUME = {9},
 PUBLISHER = {Akademie Verlag, Berlin},
      YEAR = {1996},
     PAGES = {325},
}

\bib{Gwilliam-Grady}{article}{
        AUTHOR = {Grady, R.},
	AUTHOR = {Gwilliam, O.},
	TITLE = {One-dimensional {C}hern-{S}imons theory and the {$\hat A$} genus},
	JOURNAL = {Algebr. Geom. Topol.},
	VOLUME = {14},
	YEAR = {2014},
	NUMBER = {4},
	PAGES = {2299--2377},
}


\bib{GLL}{article}{
    AUTHOR = {Grady, R.},
    author = {Li, Q.},
    author = {Li, S.},
     TITLE = {Batalin-{V}ilkovisky quantization and the algebraic index},
   JOURNAL = {Adv. Math.},
    VOLUME = {317},
      YEAR = {2017},
     PAGES = {575--639},
}




\bib{Kapranov}{article}{
    AUTHOR = {Kapranov, M.},
     TITLE = {Rozansky-{W}itten invariants via {A}tiyah classes},
   JOURNAL = {Compositio Math.},
    VOLUME = {115},
      YEAR = {1999},
    NUMBER = {1},
     PAGES = {71--113},
}



\bib{Karabegov96}{article}{
    AUTHOR = {Karabegov, A.V.},
     TITLE = {Deformation quantizations with separation of variables on a
              {K}\"{a}hler manifold},
   JOURNAL = {Comm. Math. Phys.},
    VOLUME = {180},
      YEAR = {1996},
    NUMBER = {3},
     PAGES = {745--755},
}

\bib{Karabegov00}{incollection}{
    AUTHOR = {Karabegov, A.V.},
     TITLE = {On {F}edosov's approach to deformation quantization with
              separation of variables},
 BOOKTITLE = {Conf\'{e}rence {M}osh\'{e} {F}lato 1999, {V}ol. {II} ({D}ijon)},
    SERIES = {Math. Phys. Stud.},
    VOLUME = {22},
     PAGES = {167--176},
 PUBLISHER = {Kluwer Acad. Publ., Dordrecht},
      YEAR = {2000},
}

\bib{Karabegov07}{article}{
    AUTHOR = {Karabegov, A.V.},
     TITLE = {A formal model of {B}erezin-{T}oeplitz quantization},
   JOURNAL = {Comm. Math. Phys.},
    VOLUME = {274},
      YEAR = {2007},
    NUMBER = {3},
     PAGES = {659--689},
}

\bib{Karabegov}{article}{
    AUTHOR = {Karabegov, A.V.},
    author = {Schlichenmaier, M.},
     TITLE = {Identification of {B}erezin-{T}oeplitz deformation
              quantization},
   JOURNAL = {J. Reine Angew. Math.},
    VOLUME = {540},
      YEAR = {2001},
     PAGES = {49--76},
}



\bib{Kontsevich2}{article}{
    AUTHOR = {Kontsevich, M.},
     TITLE = {Rozansky-{W}itten invariants via formal geometry},
   JOURNAL = {Compositio Math.},
    VOLUME = {115},
      YEAR = {1999},
    NUMBER = {1},
     PAGES = {115--127},
}


\bib{Lu-Shiffman}{article}{
    AUTHOR = {Lu, Z.},
    author = {Shiffman, B.},
     TITLE = {Asymptotic expansion of the off-diagonal {B}ergman kernel on
              compact {K}\"{a}hler manifolds},
   JOURNAL = {J. Geom. Anal.},
    VOLUME = {25},
      YEAR = {2015},
    NUMBER = {2},
     PAGES = {761--782},
}

\bib{Ma-Ma}{article}{
    AUTHOR = {Ma, X.},
    author = {Marinescu, G.},
     TITLE = {Berezin-{T}oeplitz quantization on {K}\"{a}hler manifolds},
   JOURNAL = {J. Reine Angew. Math.},
    VOLUME = {662},
      YEAR = {2012},
     PAGES = {1--56},

}

\bib{Melrose}{article}{
    AUTHOR = {Melrose, R.},
     TITLE = {Star products and local line bundles},
   JOURNAL = {Ann. Inst. Fourier (Grenoble)},
    VOLUME = {54},
      YEAR = {2004},
    NUMBER = {5},
     PAGES = {1581--1600, xvi, xxii},
}

\bib{Nest-Tsygan}{article}{
    AUTHOR = {Nest, R.},
    author = {Tsygan, B.},
     TITLE = {Remarks on modules over deformation quantization algebras},
   JOURNAL = {Mosc. Math. J.},
    VOLUME = {4},
      YEAR = {2004},
    NUMBER = {4},
     PAGES = {911--940, 982},
}

\bib{Neumaier}{article}{
    AUTHOR = {Neumaier, N.},
     TITLE = {Universality of {F}edosov's construction for star products of
              {W}ick type on pseudo-{K}\"{a}hler manifolds},
   JOURNAL = {Rep. Math. Phys.},
    VOLUME = {52},
      YEAR = {2003},
    NUMBER = {1},
     PAGES = {43--80},
}

\bib{Schlichenmaier}{article}{
	AUTHOR = {Schlichenmaier, M.},
	TITLE = {Deformation quantization of compact K\"ahler manifolds by Berezin-Toeplitz quantization},
	BOOKTITLE = {Conf\'erence Mosh\'e Flato 1999, Vol. II (Dijon)},
	SERIES = {Math. Phys. Stud.},
	VOLUME = {22},
	PAGES = {289–306},
	PUBLISHER = {Kluwer Acad. Publ., Dordrecht},
	YEAR = {2000},
}




\bib{Tian}{article}{
    AUTHOR = {Tian, G.},
     TITLE = {On a set of polarized {K}\"{a}hler metrics on algebraic manifolds},
   JOURNAL = {J. Differential Geom.},
    VOLUME = {32},
      YEAR = {1990},
    NUMBER = {1},
     PAGES = {99--130},
}

\bib{Tsygan}{article}{
    AUTHOR = {Tsygan, B.},
     TITLE = {Oscillatory modules},
   JOURNAL = {Lett. Math. Phys.},
    VOLUME = {88},
      YEAR = {2009},
    NUMBER = {1-3},
     PAGES = {343--369},
}




\end{biblist}
\end{bibdiv}
\end{document}